\documentclass[12pt]{article}
\usepackage{amssymb}
\usepackage{graphicx}
\usepackage{breqn}
\usepackage{latexsym}
\def\part#1{\frac{\partial\phantom{q}}{\partial#1}}

\newenvironment{rmk}{\begin{trivlist}\item[]{\bf Remark:} }
{\end{trivlist}}
\newenvironment{rmks}{\begin{trivlist}\item[]{\bf Remarks:} }
{\end{trivlist}}
\newenvironment{ex}{\begin{trivlist}\item[]{\bf Example:} }
{\end{trivlist}}
\newenvironment{prf}{\begin{trivlist}\item[]{\bf Proof:} }
{\hfill $\Box$ \end{trivlist}}

\newtheorem{thm}{Theorem}
\newtheorem{definition}{Definition}
\newtheorem{prp}[thm]{Proposition}
\newtheorem{lemma}[thm]{Lemma}

\newcommand{\lie}[1]{\mathfrak{#1}}
\def\End{\mathop{\rm End}\nolimits}
\def\Prym{\mathop{\rm Prym}\nolimits}
\def\Hom{\mathop{\rm Hom}\nolimits}

\def\ker{\mathop{\rm ker}\nolimits}
\def\coker{\mathop{\rm coker}\nolimits}

\def\Pic{\mathop{\rm Pic}\nolimits}

\def\deg{\mathop{\rm deg}\nolimits}

\def\tr{\mathop{\rm tr}\nolimits}

\def\ad{\mathop{\rm ad}\nolimits}
\def\Res{\mathop{\rm Res}\nolimits}

\def\Nm{\mathop{\rm Nm}\nolimits}

\def\ch{\mathrm {ch}}

\newcommand{\R}{\mathbf{R}}
\newcommand{\C}{\mathbf{C}}

\newcommand{\Z}{\mathbf{Z}}

\newcommand{\PP}{{\rm P}}

\begin{document}
\title{Critical loci for Higgs bundles}
\author{Nigel Hitchin\\[5pt]}

\maketitle

\section{Introduction}
 The moduli space of Higgs bundles over a curve has played a fundamental role in recent work relating Langlands duality and mirror symmetry \cite{KW}. Through its structure as an integrable system the duality of the abelian varieties which are the generic fibres implements the symmetry. Less work has been done on the singular fibres and this paper attempts to shed more light on these.
Our main focus is on the {\it critical locus}, the points of which form part of the singular fibres. Under certain genericity conditions this  is a subintegrable system with its own family of abelian varieties as fibres.

Recall that an {\it integrable system} is  a symplectic manifold $M^{2m}$ with a proper map $F:M\rightarrow B$ to a manifold $B^m$, the base of the system, which over an open set in $B$ is a  fibration by Lagrangian submanifolds. Properness implies that a generic fibre is a torus. Local coordinates on $B$ provide functions $f_1,\dots, f_m$ on $M$ which Poisson commute. 

The subspace $C_d\subset M$ at which the rank of the derivative of $F$ drops rank from $m$ to $m-d$ we call the $d$th critical locus. For an analytic integrable system there is a notion of nondegeneracy  of a  point on the critical locus. In the literature  of integrable systems (see for example \cite{Vey},\cite{Ito},\cite{Bol}) it is a theorem that the nondegenerate points of $C_d$  form a symplectic submanifold of codimension $2d$  which inherits the structure of a sub-integrable system.

In this paper we shall apply this criterion and study the consequences in the case of the moduli space of rank $2$ Higgs bundles over a curve $\Sigma$ of genus $g \ge 2$. This is a holomorphic integrable system and so analyticity automatically holds. The symplectic manifold is the moduli space ${\mathcal M}$ of stable pairs $(V,\Phi)$ where $V$ is a rank 2 bundle of fixed determinant and $\Phi$ a holomorphic section of $\End_0V\otimes K$.

 We show that the critical locus occurs where the Higgs field vanishes on a divisor $D$ of degree $d$, and the nondegenerate points correspond to $\Phi$ having simple zeros and semisimple derivative at those zeros.  
It then follows that we can identify the torus fibres of the subintegrable system as Prym varieties of spectral curves for the $K(-D)$-twisted Higgs bundle $(V,\Phi/s)$ where $s$  is the canonical section of ${\mathcal O}(D)$ vanishing on $D$. The degree $d$ of $D$ varies between $1$ and $2g-2$.

The key technique is to recognize that the Poisson commuting functions are given by Dolbeault cohomology classes in $H^1(\Sigma, K^{-1})$ which for a critical point have representatives which are supported in a neighbourhood of $D$. Then computations of Hessians can be localized and related to the local behaviour of $\Phi$ around its zeros.

Nondegeneracy means that the spectral curve, which defines a critical {\it value} of $F$,  has only ordinary double points. The critical locus intersects only part of the fibre over the critical value and we study for $d=1$  the rest of the fibre and its relationship with the critical locus, in particular  using Hecke curves, which fit naturally into the context of Higgs fields with a zero.  

Having identified the structure of  the critical locus, the second part of the paper studies some topological and differential-geometric issues. In particular we show that the critical locus ${\mathcal C}_1$ is an analytic cycle Poincar\'e dual to the second Chern class of ${\mathcal M}$ and also determine the multiplicities of the components of the nilpotent cone.

Since ${\mathcal M}$ has a natural hyperk\"ahler metric it is an obvious question to ask if the induced metric on the symplectic subvariety ${\mathcal C}_d$ is hyperk\"ahler. We show that for $d=1$ this is not so, but in fact at the other extreme $d=2g-2$ the induced metric is indeed a flat hyperk\"ahler metric. This has some bearing on the issue of BBB-branes in this moduli space, since these are expected to be supported on hyperk\"ahler submanifolds.

Finally we study the case of genus $2$  by identifying the image of the abelian varieties  for $d=1$ and $2$ in the moduli space of semistable bundles which is in this case the projective space $\PP^3$. For $d=1$ the tori are described by a pencil of singular quartic surfaces, Kummer surfaces of Jacobians, and for $d=2=2g-2$ by lines whose double cover is an elliptic curve. Using  the explicit analytical expressions for the integrable system obtained in  \cite{Lor} we have equations for these.

The author wishes to thank A. Oliveira for useful conversations and EPSRC and ICMAT for support during the preparation of this work.

\section{A classical example}\label{classic}
As a concrete guide to what to expect in an integrable system consider the classical case of the geodesics on an ellipsoid 
$$\frac{x^2}{a^2}+\frac{y^2}{b^2}+\frac{z^2}{c^2}=1$$
where $0<a<b<c$. Here the symplectic manifold is the cotangent bundle $T^*S^2$. The   Riemannian metric defines a function $f_1$ on cotangent vectors, and there is a second quadratic function $f_2$  (discovered by Bessel) which Poisson commutes with the first. The geodesic flow is the Hamiltonian vector field $X_1$ of $f_1$. Then $F=(f_1,f_2)$ gives a map to $\R^2$ which is the integrable system. The generic fibre is a 2-dimensional torus. 

 Umbilical points on a surface are where the trace-free component of the second fundamental form vanishes. This is a smooth section of the complex line bundle $K^{-2}$ in general and on the ellipsoid $\cong S^2$ this is of degree $4$. The four umbilical points  consist of  two pairs of  points $\pm {u},\pm{v}\in \R^3$ on the ellipse $y=0$. A singular fibre of the integrable system is defined by $(f_2=0, f_1 = const)$ and (see \cite{CdV}) consists of the orbits of the geodesic flow  on $T^*S^2$ which project to geodesics passing through the umbilical points. These are closed geodesics and are given by two tori (one for the geodesics through $\pm { u}$, one for $\pm { v}$) and they intersect in two circles corresponding to the closed geodesic $y=0$ joining  $ {u}$ and ${ v}$ with the two orientations. 
 
 The  critical locus in this fibre is the pair of circles and, as $f_1$ varies (i.e. we consider the geodesic flow on circle bundles of different size) we obtain a real symplectic submanifold which is a  bundle over the line with fibre a pair of circles. The point to notice here is that the singular fibre  is the image of a smooth manifold -- a disjoint pair of tori -- and the singular locus of the image  is the critical locus. This is itself  a smooth subintegrable system (away from $f_1=0=f_2$).
 \section{The Higgs bundle system}
 \subsection{The critical locus}
We restrict ourselves for convenience to $SL(2,\C)$-Higgs bundles over a compact Riemann surface $\Sigma$ of genus $g\ge 2$, i.e. pairs $(V,\Phi)$ such that $V$ is a rank 2 holomorphic vector bundle  with $\Lambda^2V\cong {\mathcal O}$ and $\Phi\in H^0(\Sigma,\End_0V\otimes K)$ with $\tr \Phi=0$.  The moduli space ${\mathcal M}$ of stable Higgs bundles of this type \cite{Hit1},\cite{Sim} has dimension $6g-6$. The map $F$ defined by $(V,\Phi)\mapsto \tr\Phi^2/2$ gives the integrable system, the base ${\mathcal B}=H^0(\Sigma, K^2)$ being in this case a vector space. 

We have the characteristic polynomial $\det(x-\Phi)=x^2-q$ where $q=\tr\Phi^2/2$ and then $q\in {\mathcal B}$ defines a curve $x^2-q=0$, the spectral curve $S$. This lies in the total space of the canonical bundle $K$ of $\Sigma$  and the projection $\pi:K\rightarrow \Sigma$ restricts to $S$ to define a double covering. We shall only be considering the case where $S$ is irreducible which means in particular that  there are no invariant subbundles. This in turn implies that $(V,\Phi)$ is stable and we are at a smooth point of ${\mathcal M}$. 

When $S$ is nonsingular, the Lagrangian fibre of $F:{\mathcal M}\rightarrow {\mathcal B}$ is isomorphic to the Prym variety of $S\rightarrow \Sigma$, the kernel of the norm map $\Nm:\Pic^0(S)\rightarrow\Pic^0(\Sigma)$. Equivalently, using the involution  $\sigma(x)=-x$, the Prym variety consists of the anti-invariant part of $\Pic^0(S)$, line bundles $L$ such that $\sigma^*L\cong L^{-1}$. The Higgs bundle is recovered from  $L$  by $V=\pi_*(L\otimes \pi^*K^{1/2})$ and $\Phi=\pi_*x$ where $x:L\rightarrow L\otimes \pi^*K$ is the action on $L$ of the tautological section $x$ of $\pi^*K$ over  $K$.

The critical locus consists of the points at which the derivative of $F$ drops rank, so we need a geometrical way of understanding this derivative.

From  \cite{Bis} we know that at a smooth point $(V,\Phi)$ the tangent space to ${\mathcal M}$ can be understood as  the first hypercohomology group of the sequence of sheaves 
               \begin{equation}
               {\mathcal O}(\End_0V)\stackrel {[\Phi , -]}\rightarrow {\mathcal O}(\End_0V\otimes K).
               \label{complex}
               \end{equation}
               
               There are two spectral sequences associated to the hypercohomology of a double complex. The most obvious one as in \cite{Hit1} gives 
                $$\rightarrow H^0(\Sigma, \End_0V\otimes K)\rightarrow {\mathbb H}^1\rightarrow H^1(\Sigma, \End_0V)\rightarrow $$
                which, when $V$ is stable, is just the structure of the tangent space to the cotangent bundle of the moduli space ${\mathcal N}$ of stable bundles.
                
               Here instead we take the one whose $E^2$ term involves the cohomology sheaves of the complex above: namely the kernel and cokernel sheaves in (\ref{complex}). By stability  ${\mathbb H}^0$ and  ${\mathbb H}^2$ vanish and we obtain an exact sequence 
         \begin{equation}
               0\rightarrow H^1(\Sigma, \ker \Phi)\rightarrow {\mathbb H}^1\rightarrow H^0(\Sigma, \coker \Phi)\rightarrow 0.
                   \label{hyper}
               \end{equation}         
               
               \begin{prp} \label{nonvan} If $\Phi$ is non-vanishing everywhere then $(V,\Phi)$ is not on the critical locus.
                \end{prp}
                \begin{prf}
         Consider first the kernel, the sheaf of centralizers of $\Phi$.  In the Lie algebra $\lie{sl}(2,\C)$ an element $b$ commutes with a nonzero element $a$ if and only if $b=\lambda a$. 
         So if $\Phi$ is everywhere non-zero then $\Phi:K^{-1}\rightarrow \End_0V$ identifies the sheaf $\ker \Phi$ with ${\mathcal O}(K^{-1})$. 
          The image of $\ad a$ is the orthogonal complement of $a$ using the invariant  inner product $\tr ab$  so the cokernel in $\lie{sl}(2,\C)$ is the dual of the kernel. In $\End_0V\otimes K$ this makes the cokernel ${\mathcal O}(K^2)$. Then the exact sequence above expresses the tangent space as 
 $$0\rightarrow H^1(\Sigma, K^{-1})\rightarrow T_{(V,\Phi)}\rightarrow H^0(\Sigma, K^2)\rightarrow 0.$$
   Concretely, the second homomorphism in (\ref{hyper}) can be seen as follows. We use Dolbeault formalism to write the holomorphic structure on $V$ as a $\bar\partial$-operator $\bar\partial_A$ and then a tangent vector to ${\mathcal M}$ is represented by $(\dot A,\dot\Phi)\in \Omega^{01}(\Sigma, \End_0V)\oplus \Omega^{10}(\Sigma, \End_0V)$ such that 
                 \begin{equation}
               \bar\partial_A\dot\Phi+[\dot A,\Phi]=0
                 \label{tangt}
               \end{equation}
  If $p:\End_0 V\rightarrow  \coker \Phi$ is the projection to the cokernel, then $p(\dot \Phi)$ is holomorphic from (\ref{tangt}) since $p([\dot A,\Phi])=0$.                Then   $\tr(\Phi\dot\Phi)=\tr(\Phi p(\dot\Phi))$ identifies the image in  $\coker \Phi\cong H^0(\Sigma, K^2)$. But this  is also clearly the derivative of $F=\tr\Phi^2/2$. By exactness in (\ref{hyper}) there are no critical points if $\Phi$ is non-vanishing.
  \end{prf}
  Note that  in this situation the kernel of $DF$ is isomorphic to $H^1(\Sigma,K^{-1})$. Pull back a class $a$ to $\pi^*a\in H^1(S,\pi^*K^{-1})$ and multiply by $x\in H^0(S,\pi^*K)$ and then $xa\in H^1(S,{\mathcal O})$ is anti-invariant by the involution $\sigma(x)=-x$. In this way the kernel of $DF$ is recognizable as  the tangent space of the Prym variety.

 \begin{prp} \label{dim} If $\Phi$ vanishes on a divisor $D$, then $1\le \deg D\le 2g-2$.
 \end{prp}
 \begin{prf}
         Suppose  the zeros of $\Phi$ are given by the divisor $D=m_1x_1+\cdots + m_{\ell}x_{\ell}$ then           
          locally around a zero $x_i\in \Sigma$ we have $\Phi=z^m\phi$ where $\phi$ is non-vanishing and then  $\det \Phi=z^{2m}\det\phi$ and so when $\phi(0)$ is semisimple we have a zero of the quadratic differential $q=\det\Phi$ of order $2m$ and in the nilpotent case $\ge (2m+1)$. Suppose $x_k$ for $1\le k\le m$ are the nilpotent zeros and $x_k$ for $m+1\le k\le \ell$ the semisimple ones. Then since $\deg K^2=4g-4$,
          $$4g-4\ge \sum_{k=1}^{m}(2m_k+1)+\sum_{k=m+1}^{\ell}2m_k=2\deg D+m$$
          and so $\deg D\le  2g-2$. 
      \end{prf}           
                   
Therefore  at a point on the critical locus $\Phi$ must vanish on a divisor $D$ and so $\Phi=s\phi$ where $s$ is the canonical section of ${\mathcal O}(D)$ vanishing on $D$ and $\phi$ is a non-vanishing section of $\End_0V\otimes K(-D)$.

\subsection{Nondegeneracy}\label{nondegen}

 If the Hessian of a function at a critical point is nondegenerate in the usual sense, then that critical point  is isolated. But  for an integrable system,  because of the commuting vector fields, this can only occur  when it is critical for all functions.  
 There is however a notion of nondegeneracy of critical points for integrable systems \cite{Vey}, \cite{Ito}, \cite{Bol}. If $(M^{2m},\omega)$ is the symplectic manifold we consider the induced action on the tangent space $T_x$ at $x\in M$ of the  vector fields $X_1,\dots, X_d$ that vanish at $x$. They act as a commutative subalgebra of $\End T_x$   preserving the restriction of the symplectic form. There is a distinguished subspace $E\subset T_x$ of dimension $(m-d)$ spanned by {\it all} the Hamiltonian vector fields  and this is  preserved  since the vector fields  commute. Since $\omega(X_i,X_j)=0$ for all Hamiltonian vector fields, $E$ is isotropic and hence $E^{\perp}/E$ is a vector space of dimension $(2m-(m-d))-(m-d)=2d$ which has an induced symplectic form. Then 
   $X_1,\dots, X_d$ span a commutative subalgebra of $\lie{sp}(2d,\C)$. 
   
   \begin{definition} \label{defn}The point $x$ is called {\it nondegenerate} if the induced action of $X_1,\dots, X_d$ spans  a Cartan subalgebra of  $\lie{sp}(2d,\C)$
   \end{definition}
  If $F=(f_1,\dots,f_m)$ defines the integrable system then identifying the Lie algebra $\lie{sp}(2d,\C)$ with quadratic functions on $\C^{2d}$ this subalgebra is spanned by the Hessians of the Hamiltonian functions $f_1,\dots,f_d$ for which the point is critical.  This is essentially the induced action at a fixed point on the local symplectic quotient by the action of the vector fields $X_{d+1},\dots, X_m$.
  
  \begin{ex} If $M=\C^2$ with symplectic form $dz_1\wedge dz_2$ then the function $z_1z_2$ has a non-degenerate critical point at the origin. It is the moment map for the semisimple action of $\C^*$ given by $(z_1,z_2)\mapsto (\lambda z_1,\lambda^{-1}z_2)$. The unipotent action of $\C$ given by $(z_1,z_2)\mapsto (z_1+\lambda z_2,z_2)$ has moment map $z_2^2/2$ and the critical locus of this is the line $z_2=0$, a degenerate locus.
    \end{ex}
    
    To reformulate this criterion slightly, let  $f_1,\dots,f_{r}$  be components of the function $F$ such that $df_i$ are linearly independent at $x$ and $X_1,\dots, X_{r}$ be the corresponding Hamiltonian vector fields so that they span $E$. The condition $\omega(X,X_i)=0$   is equivalent to $df_i(X)=0$ for all $1\le i\le r$ hence   $E^{\perp}= \ker DF_x$.   The symplectic form restricted to  $E^{\perp}$ will have  $E$  as its degeneracy subspace. We next calculate $E^{\perp}/E$ for Higgs bundles.
     
     Recall that the symplectic form on ${\mathcal M}$ is defined on representatives by 
     $$\omega((\dot A_1,\dot\Phi_1),(\dot A_2,\dot\Phi_2))=\int_{\Sigma}\tr(\dot A_1\dot\Phi_2)-\int_{\Sigma}\tr(\dot A_2\dot\Phi_1).$$
     
     Let $(\dot A,\dot\Phi)$ represent a tangent vector in $\ker DF$. Then $\tr \Phi\dot\Phi=0$.  But  $\phi=\Phi/s$ is smooth and non-vanishing so $\dot \Phi=[\phi,\psi]$ for a smooth section $\psi$ of $\End_0V (D)$, well-defined up to a multiple of $\phi$.  Now $\bar\partial \dot\Phi+[\dot A,\Phi]=0$ yields
      \begin{equation}
    [\Phi,\dot A]=\bar\partial \dot \Phi=[\Phi,\bar\partial \psi/s]
    \label{psi}
    \end{equation}
     and so $\dot A=\bar\partial \psi/s+\lambda \Phi/s$. 
           Consider then, using $\tr \Phi\dot\Phi_2=0$
    $$\int_{\Sigma}\tr(\dot A_1\dot\Phi_2)=\int_{\Sigma}\tr\left(\frac{\bar\partial \psi_1\dot\Phi_2}{s}+\lambda_1 \frac{\Phi\dot\Phi_2}{s}\right)=\int_{\Sigma}\tr\left(\frac{\bar\partial \psi_1[\Phi,\psi_2]}{s^2}\right)$$
  and  
    $$\omega((\dot A_1,\dot\Phi_1),(\dot A_2,\dot\Phi_2))=\int_{\Sigma}\tr\left(\frac{\bar\partial \psi_1[\Phi,\psi_2]-\bar\partial \psi_2[\Phi,\psi_1]}{s^2}\right)=\int_{\Sigma}\bar\partial\left(\frac{\tr( \psi_1[\Phi,\psi_2])}{s^2}\right).$$
    If $s$ were non-vanishing, so that $\ker DF$ is the tangent space to a  regular fibre, then by Stokes' theorem the integral would vanish. This would simply confirm  that the fibre is Lagrangian. Near a  simple zero $z=0$ of $s$ we have from above that $\Phi=z\phi$ and  then the integral has a contribution 
    \begin{equation}
    2\pi i\tr(\psi_1[\phi,\psi_2])(0).
    \label{simp}
    \end{equation}
   The expression (\ref{simp}) is a multiple of the canonical  symplectic form  $\omega$ on the coadjoint orbit of $\phi(0)$ in ${\lie sl}(2,\C)$ and so the symplectic vector space $E^{\perp}/E$, which  is  $\ker DF$ modulo the degeneracy subspace, is isomorphic to the  direct sum of $\deg D$ two-dimensional tangent spaces of a  coadjoint orbit of $SL(2,\C)$. Then $2\deg D=\dim  E^{\perp}/E = \dim \ker DF-(6g-6-\dim\ker DF)$ hence $\dim \ker DF=3g-3+\deg D$.
   
   Near a higher order zero  we have $\Phi=z^k\phi$ where $\phi$ is non-vanishing and from $\dot A=\bar\partial \psi/s+  \lambda \Phi$ we see that $\psi$ is holomorphic modulo $z^k$  so we can write the local formula for the symplectic form as  
  $$\frac{1}{2\pi i}\Res_{z=0}\frac{\tr(\psi_1[\phi,\psi_2])}{z^k}$$
  where $\phi,\psi_1,\psi_2$ are polynomials of degree $(k-1)$ with values in $\lie{sl}(2,\C)$. This is the tangent space of an orbit in a truncated loop algebra. It is a symplectic vector space of dimension  $2k$. 
  
 An immediate consequence is  the following:
 \begin{prp} \label{image} Let $(V,\Phi)$ be a Higgs bundle such that $\Phi$ vanishes on a divisor $D$. Then $(V,\Phi)$ lies on the $d$th critical locus where $d=\deg D$ and the image of $DF$  is $H^0(\Sigma, K^2(-D))\subset H^0(\Sigma, K^2)$. 
   \end{prp}
  \begin{prf} The dimension of $\ker DF$ is calculated as above to be $3g-3+\deg D$ so the rank is $3g-3-\deg D$ and $(V,\Phi)$ is on the critical locus for $d=\deg D$.
  
   If $\Phi$ vanishes on $D$ then clearly $\tr \Phi\dot\Phi\in H^0(\Sigma, K^2(-D))$ so the image is a $(3g-3-d)$-dimensional subspace. But from Proposition \ref{dim} $\deg K^2(-D)\ge 2g-2$ and so by Riemann-Roch $\dim H^0(\Sigma, K^2(-D)) =3g-3-d$ unless $D$ is a canonical divisor. In the latter case $\phi\in H^0(\Sigma, \End_0V)$ and an eigenspace defines an invariant degree zero subbundle which contradicts stability. It follows that $DF$ surjects onto $H^0(\Sigma, K^2(-D))$.
  \end{prf}

   \subsection{Critical Hamiltonians}
    For the  Higgs bundle integrable system  $H^1(\Sigma, K^{-1})$ is the dual space of the base ${\mathcal B}=H^0(\Sigma, K^2)$ by Serre duality, so these $3g-3$ linear functions are Hamiltonian functions $f$. If $f$ is given by a class $a\in H^1(\Sigma, K^{-1})$ represented by  $\beta\in \Omega^{01}(\Sigma, K^{-1})$ then 
    $$f=\frac{1}{2}\int_{\Sigma} \beta\tr\Phi^2.$$
    At a critical point at which $\Phi$ vanishes on a divisor $D$  we need to know the $d$-dimensional subspace of the space of linear functions $H^1(\Sigma, K^{-1})$ which have this as a critical point. 
   
 We saw in Proposition \ref{image} that the image of $DF$ is $H^0(\Sigma,K^2(-D))$ embedded in $H^0(\Sigma, K^2)$ by the product with $s$. By Serre duality its annihilator is the kernel of $H^1(\Sigma, K^{-1})\stackrel{s}\rightarrow H^1(\Sigma , K^{-1}(D))$. Since $d=\deg D \le 2g-2$ and $D$ is not a canonical divisor, $H^0(\Sigma, K^{-1}(D))=0$ and  we have the exact sequence 
 $$0\rightarrow H^0(D, K^{-1}(D))\stackrel{\delta}\rightarrow H^1(\Sigma, K^{-1})\stackrel{s}\rightarrow H^1(\Sigma , K^{-1}(D))\rightarrow 0$$
 and so the functions $f$ are given by the $d$-dimensional 
 image of 
$H^0(D, K^{-1}(D))$.

We can represent such classes in $H^1(\Sigma, K^{-1})$ by forms which are localized around $D$ and this will give in particular expressions for the Hessian of these functions.  Given a section in $H^0(D, K^{-1}(D))$ a Dolbeault representative of the class in $H^1(\Sigma, K^{-1})$ can be taken of the form $\bar\partial \alpha/s$ where $\alpha$ is supported in a  neighbourhood of the zeros of the section $s$. We extend the section holomorphically from $D$ to a neighbourhood and then using a bump function construct a $C^{\infty}$ section $\alpha$ of the line bundle $K^{-1}(D)$.
 Then the class is represented by $\bar\partial\alpha/s$, a  sum of forms which have local expressions 
          $$\frac{1}{z^{m_k}}\bar\partial \alpha_k$$ 
          around a zero $z=0$ of $\Phi$ of order $m_k$.

          To calculate the Hessian  for the class $a=[\bar\partial \alpha/s]\in H^1(\Sigma,K^{-1})$ we  consider a second order deformation  $$A(t)=A+t\dot A+ t^2\ddot A, \qquad\Phi(t)=\Phi+t\dot \Phi+ t^2\ddot \Phi$$  so that
        \begin{equation}
        \bar\partial_A\dot\Phi+[\dot A,\Phi]=0\qquad \bar\partial_A\ddot\Phi+[\dot A,\dot\Phi]+[\ddot A,\Phi]=0.
        \label{vary}
        \end{equation}
    Then, 
    $$\int_{\Sigma}\tr \Phi(t)^2  \bar\partial \alpha  /s=\int_{\Sigma}(\tr \Phi^2 +2t\tr \Phi\dot\Phi+t^2(\tr \dot\Phi^2+2\Phi\ddot\Phi)) \bar\partial \alpha  /s+o(t^3)$$
    and the second variation at a critical point is the term
    $$\int_{\Sigma}(\tr \dot\Phi^2+2\tr \Phi\ddot\Phi) \bar\partial \alpha  /s.$$
    Since $\Phi$ vanishes on $D$, $\Phi\alpha/s$ is smooth and  Stokes' theorem gives 
 \begin{equation}
 \int_{\Sigma}\tr \dot\Phi^2 \bar\partial \alpha /s-2\int_{\Sigma}\tr(\Phi  \bar\partial_A \ddot\Phi)\alpha /s.
 \label{f1}
 \end{equation}
    But from (\ref{vary}) 
    $$\tr(\Phi  \bar\partial_A \ddot\Phi)=-\tr(\Phi ([\dot A,\dot\Phi]+[\ddot A,\Phi]))= -\tr(\Phi [\dot A,\dot\Phi])=\tr([\dot A,\Phi] \dot\Phi)=-\tr(\bar\partial_A\dot \Phi\dot\Phi)$$
    and we can rewrite (\ref{f1}) as 
    $$ \int_{\Sigma}\tr \dot\Phi^2 \bar\partial \alpha /s+\int_{\Sigma}\bar\partial (\tr \dot\Phi^2)\alpha /s$$
    and by Stokes' theorem this is a boundary integral.  
    Since $\alpha$ is holomorphic in a small neighbourhood of the zeros of $s$ the integral  is equal to
    \begin{equation}
    2\pi i \sum_{s(x_k)=0} \Res_{x=x_k} \frac{\alpha\tr \dot\Phi^2}{s}.
    \label{hess}
    \end{equation}
    In particular the contribution from a simple zero is a multiple  of $\tr \dot\Phi^2$.
    \begin{rmk}
    Note that this expression is gauge-invariant -- replacing $\dot\Phi$ by $\dot\Phi+[\psi,\Phi]$ gives $\tr (\dot\Phi+[\psi,\Phi])^2=\tr\dot\Phi^2+2\tr( \dot\Phi[\psi,\Phi])+\tr  [\psi,\Phi]^2$ and the last two terms are divisible by $s$.
    \end{rmk}

 Using this localized calculation we can determine nondegeneracy: 
  
  \begin{prp}  \label{non} Let $(V,\Phi)$ be a critical point  of $F$ on the moduli space of $SL(2,\C)$-Higgs bundles. Then it is nondegenerate in the above sense if and only if $\Phi$ has simple zeros on the divisor $D$ and its derivative is semisimple at those zeros. 
     \end{prp}
  \begin{prf} In the previous section we identified the symplectic vector space $E^{\perp}/E$ with a direct sum of tangent spaces to coadjoint orbits. The action of the commuting vector fields on this space is given by the Hessians of the Hamiltonian functions, interpreted as quadratic functions on $E^{\perp}/E$, elements in the Lie algebra of symplectic   transformations. 
  
  At a simple zero the symplectic form was 
 $$ 2\pi i\tr(\psi_1[\phi,\psi_2])(0).$$

   The Hamiltonian vector field on the coadjoint orbit of $a\in \lie{sl}(2,\C)$ defined by $a$ itself using the Killing form has a zero at the point $a$ on this orbit and its Hessian is the symmetric form 
  $H(u,v)=\omega([a, u],v)$. In our case this is, up to a multiple, 
  $$\tr([\phi,\psi_1][\phi,\psi_2])=\tr (\dot \Phi_1\dot\Phi_2)$$
  which is the Hessian as in (\ref{hess}). Thus when $\Phi$ has simple zeros and semisimple derivative the action on $E^{\perp}/E$ is via $d$ commuting semisimple elements and hence a Cartan subalgebra. This is the nondegeneracy definition. If $\phi(0)$ is nilpotent then it does not form part of a Cartan subalgebra.
  
 For a zero of order $k>1$  there are critical functions for which $\alpha=z^i$ with $0<i\le k-1$ and their 
Hessians  are given by 
  $$\Res_{z=0} \frac{z^i\tr \dot\Phi_1\dot\Phi_2}{z^k}=\Res_{z=0} \frac{z^i\tr( [\phi,\psi_1][\phi,\psi_2])}{z^k}$$
  where we work with a Lie algebra of polynomials modulo $z^k$. 
  For $i>0$, $(z^i\phi)^k=0$ modulo $z^k$ which means that the  Hessian for $z^i$ is nilpotent in its action. 
  
  Consequently only when $k=1$ and $\phi(0)$ is semisimple  do we get a Cartan subalgebra. 

   \end{prf}
   
   \begin{rmk} Another interpretation of Proposition \ref{non} is that nondegenerate critical points are those in a fibre for which  the quadratic differential $q=\det \Phi$ has  zeros of multiplicity at most two, or equivalently the spectral curve $S$ defined by $x^2-q=0$ has only ordinary double points, a natural notion of genericity.
   \end{rmk}
   \subsection{${\mathcal C}^0_d$ as an integrable system} \label{Cint}
   
   Let ${\mathcal C}_d$ denote the $d$th critical locus and  ${\mathcal C}^0_d$ the open submanifold of nondegenerate points. According to the general theorem \cite{Ito}, this is an integrable subsystem of ${\mathcal M}$ of dimension $6g-6-2d$ with complex torus fibres. We shall analyze its structure  in this section.

We have seen that at if $(V,\Phi)$ is a critical point then $\Phi=\phi/s$ where $\phi\in H^0(\Sigma, K(-D))$.  Moduli spaces of Higgs bundles  twisted by a general line bundle $L$ instead of $K$ have long been constructed \cite{Nit}, the chief difference being that we no longer have a symplectic form defined on them. The function $\tr \phi^2/2$ defines   a proper map  to $H^0(\Sigma, L^2)$ as before, though this now has a different  dimension from the fibre. When $\deg L\ge 2g-2$ Riemann-Roch gives more regularity, but in our case when $\deg K(-D)< 2g-2$ we still have an open set with  spectral curve constructions.

 At a stable point the tangent space is again given by a hypercohomology space ${\mathbb H}_d^1$ where now
   \begin{equation}
   \label{subhyper}
     0\rightarrow H^1(\Sigma, K^{-1}(D))\rightarrow {\mathbb H}_d^1\rightarrow H^0(\Sigma, K^2(-2D))\rightarrow 0.
\end{equation}
   and since for us $\phi$ is now everywhere non-vanishing, this describes as in  Proposition \ref{nonvan} the derivative of $\tr\phi^2$ on a smooth fibre.
   
   The spectral curve approach is reflected by this exact sequence, and the $(3g-3-d)$-dimensional space $H^1(\Sigma, K^{-1}(D))$ is tangent to the Prym variety of the 
    spectral curve $\tilde S$ defined by $\tilde x^2-\det\phi=0$ in the total space of $K(-D)$. 
    This $(3g-3-d)$ is half the dimension of ${\mathcal C}^0_d$ to  
      show that it is Lagrangian with respect to the induced symplectic structure we only need to show it is isotropic. 
      
   But  in Section \ref{nondegen} we computed the symplectic form on the kernel of $DF$, and the fibre we are considering now has $s$ fixed, which means that  $\dot\Phi$  in Equation (\ref{psi}) is divisible by $s$. This removes the denominator in the calculation and Stokes' theorem gives zero.
   
  Nondegeneracy gives us an integrable system, and  we have just identified a smooth fibre as isomorphic to the Prym variety of $\tilde S$.  The symplectic form on ${\mathcal C}^0_d$ identifies the tangent space at a point on the base with the dual of $H^1(\Sigma, K^{-1}(D))$, the tangent space of the fibre. By Serre duality this is $H^0(\Sigma, K^2(-D))$ which we have already seen  is the image of the restriction of $DF$ to ${\mathcal C}^0_d$.
  
   The base ${\mathcal B}_d$ of this subsystem is the space of quadratic differentials with $d$ double zeros and no more. If $d\ge g-1$, Riemann-Roch no longer gives the dimension of $H^0(\Sigma, K^2(-2D))$ but its generic value is $(3g-3-2d)$, and this must hold for  the divisors coming from ${\mathcal C}^0_d$. This is the base for the $K(-D)$-twisted moduli space. The product with $s$ gives an embedding  
  $H^0(\Sigma, K^2(-2D))\subset H^0(\Sigma, K^2(-D))$ with quotient $H^0(D,K^2(-D))$ gives the infinitesimal deformations of $D$, the double zeros of the quadratic differential.
  \subsection{Calabi-Yau threefolds}
  
  The spectral curves for Higgs bundles which are critical for a specific linear function $f$ have another geometric description. 
  
  Take $f$ to be defined by  a class $a\in H^1(\Sigma, K^{-1})$ in the kernel of the map $H^1(\Sigma, K^{-1})\rightarrow H^1(\Sigma, K^{-1}(D))$ and consider the corresponding rank $2$ bundle $E$ over $\Sigma$ defined as an extension 
 $$0\rightarrow {\mathcal O}\rightarrow E\rightarrow K\rightarrow 0.$$
 Translation by the non-zero section defining the trivial subbundle is a free $\C$-action with quotient $K$, a principal $\C$-bundle defined by the class $x\pi^*a \in H^1(K,{\mathcal O})$. This total space is a $3$-manifold which is  Calabi-Yau  since the cotangent bundle is symplectic and the bundle along the $\C$-fibres is trivial. 
 
 Consider this  principal $\C$-bundle restricted to the singular spectral curve $S\subset K$. Its normalization is the smooth spectral curve $\tilde S$ for $\phi=\Phi/s$. We represent the class $a$ by a form $\bar\partial{\alpha/s}$ so the class $x\pi^*a\in H^1(K,{\mathcal O})$ has Dolbeault representative $\bar\partial{x\alpha/s}$. Near a zero of $s$ of order $m$ we have the singularity of the spectral curve given as $x^2=z^{2m}$ or $x^2=z^{2m+1}$. In the first case the normalization has a local description $x=\pm z^{2m}$ in which case $x\alpha/s$ is smooth. In the other case there is a local parametrization $x=t^{2m+1}, z=t^2$ and $x/z^m=t$ which is also smooth. Hence  the the class $a$ pulled back to  $\tilde S$ is trivial. This means the principal $\C$-bundle has a section over $\tilde S$. 
  
 It follows that each Higgs bundle in the critical locus for $f$ is represented by a line bundle over a smooth curve in a Calabi-Yau threefold. The points on this curve which are equivalent under  the $\C$-action give the singularities of $S$ in the quotient $K$.  Varying the class $a$ gives a family of such Calabi-Yau threefolds. 
 
 \begin{rmk}Note that the local structure of sheaves supported on curves in a Calabi-Yau threefold is known from DT-theory to be that of the critical locus of a function. 
 \end{rmk}

A similar situation concerning Nahm's equations is described in \cite{Hit3} where we have  $\PP^3\backslash \PP^1$ instead of the Calabi-Yau.

\subsection{The full singular fibre}\label{full}
   The critical locus is contained in the union of the singular fibres of $F:{\mathcal M}\rightarrow {\mathcal B}$, namely where the torus degenerates to a singular variety. The nondegenerate points on the critical locus form, as we have seen, a smooth submanifold, and intersect a singular fibre in a smooth torus of lower dimension. In this section we study how it is located in the full singular fibre, which gives an additional viewpoint on the approach so far. A detailed study of the singular fibres in rank $2$ is the paper \cite{Got}. 
   
   Any point in a singular fibre has associated to it a singular spectral curve $S$.  When the Higgs field $\Phi$ has no zeros then   $\coker (x-\Phi)$ on $S$ defines a line bundle even when  $S$ is singular and the  Higgs bundle is obtained  as in the smooth case by the direct image as in \cite{BNR}. We can see the Higgs field  concretely near a double zero of $\tr \Phi^2$: around a singularity $x^2=z^2$ in the spectral curve of the type we are considering a line bundle is locally trivial.   So the structure of its direct image is equivalent to  ${\mathcal O}(\pi^{-1}(U))$ as a module over ${\mathcal O}(U)$ for $U\subset \Sigma$. This consists of functions $f_0(z)+xf_1(z)$ giving $(f_0(z),f_1(z))$ as local sections of $\pi_*L=V$. The Higgs field $\Phi$ is the direct image of $x$ so $x(f_0(z)+xf_1(z))=z^2f_1(z)+xf_0(z)$ and we have 
\begin{equation}
\Phi=\pmatrix{0 & z^2\cr
                          1 & 0}dz.
                          \label{nonv}
                          \end{equation}
Note that this is non-vanishing at $z=0$. 

For simplicity assume we have just one singularity of this type on $S$ over the point $y\in \Sigma$. The singular spectral curve has a smooth normalization $\tilde S$ which is just  the spectral curve of the $K(-y)$-twisted Higgs field $\phi=\Phi/s$ and there are two points $p,q\in \tilde S$ which map to the singular point.  A line bundle on the singular curve is a line bundle $L$ on $\tilde S$ together with an identification of the  fibres over $p$ and $q$.  A change of  the identification defines a  $\C^*$-action which, as $L$ varies, generates a 
 principal $\C^*$-bundle over $\Prym(\tilde S)$. 
 
 The full fibre is compact since $F$ is proper and its compactification is provided by rank one torsion free sheaves on $S$  \cite{Sim}. When $S$ is irreducible, such a sheaf is always the direct image of a line bundle on a partial normalization, in this case $\tilde S$ \cite{Alex}. The normalization of $x^2-z^2$ around the singularity consists of two disjoint components  $x=\pm z$ with coordinate $z$ and then the direct image gives a  Higgs field vanishing at $z=0$ with the local form 
   $$\Phi=\pmatrix{z & 0\cr
                          0 & -z}dz$$
namely a simple zero at $z=0$  with   semisimple derivative. 

There are different approaches to seeing how each $\C^*$ is compactified in the moduli space by such Higgs bundles. The authors of \cite{Got} use the moduli space of generalized parabolic bundles: the identification of $L_p$ with $L_q$ has a graph  which generates a one-dimensional subspace of  $L_p\oplus L_q$ and the compactification is to the projective line $\PP(L_p\oplus L_q)$.

An alternative is to consider  the Hecke correspondence as in \cite{HwRam}. Since  the zero of the Higgs field plays a key role here we develop it further.

Let $V$ be a rank $2$ vector bundle and take a point $y\in \Sigma$ and $\alpha\in V^*_y$ a linear form on the fibre at $y$. Then define a vector bundle $V'$ by the kernel of the sheaf homomorphism:
$$0\rightarrow {\mathcal O}(V')\rightarrow {\mathcal O}(V)\stackrel{\alpha}\rightarrow {\mathcal O}_y\rightarrow 0$$
This is a new vector bundle and if $\Lambda^2V \cong {\mathcal O}(y)$ then $\Lambda^2V'$ is trivial.
In local terms, if 
$\alpha(v_1,v_2)=v_1$
then a 
local basis of $V'$ is  $e_1=(z,0),  e_2=(0,1)$. Multiples of $e_1$ form a distinguished one-dimensional subspace of $V'_y$.

Suppose now that $\Phi$ is a Higgs field for $V$, so that locally 
$$ \Phi=\pmatrix{a & b\cr
                            c & -a}$$
then in the new basis it becomes 
$$ \Phi'=\pmatrix{a & bz^{-1}\cr
                            cz & -a}$$
   and this is well-defined on $V'$ only if                          
                     $b=zb'$.
  By varying $\alpha$ it is clear that $\Phi$ lifts for all $\alpha\in V^*_y$ if $\Phi=z\phi$. In that case 
\begin{equation}
 \Phi'=\pmatrix{a'z & b'\cr
                            c'z^2 & -a'z}
                            \label{Higgsprime}
                            \end{equation}
                            where 
                            $$ \phi=\pmatrix{a' & b'\cr
                            c' & -a'}.$$
                             Hence a Higgs field on $V$ with a zero at $z=0$ defines naturally on $V'$ a Higgs field with the same characteristic polynomial.  If the characteristic polynomial is irreducible there are no invariant subbundles and so we obtain a stable Higgs bundle for all choices of $\alpha\in V^*_y$ and hence a copy of the projective line $\PP(V^*_y)$ in ${\mathcal M}$. This is a Hecke curve \cite{HwRam}. 
                             
        \begin{rmk} Since the base $B$ of the integrable system is a vector space, any projective line must lie in a fibre of $F$ and since $\PP^1$ is simply connected it cannot lie in an abelian variety hence it must be  a subvariety of a singular fibre. 
                         \end{rmk}
                         
                         There are two cases:
                         
                         \begin{itemize}
                         \item
                         If $b'(0)\ne 0$ in (\ref{Higgsprime}) then $\Phi'(0)\ne 0$ and a local change of basis will take it into the form of (\ref{nonv}). The Hecke transform  then gives a Higgs bundle which comes from a line bundle on the singular spectral curve.
                         \item
                         Suppose $b'(0)=0$, then 
                     since  $\tr \Phi^2/2= z^2({a'}^2+b'c')$, and $S$ has just an ordinary double point, we have $({a'}^{2}+b'c')(0)\ne 0$ 
                          hence $a'(0)\ne 0$.  It then follows that $\Phi'=z\phi$ with $\phi(0)$ semisimple and eigenvalues $\pm a'(0)$.This   is therefore a point on the critical locus. 
                          \end{itemize}
                          Varying  $\alpha\in V^*_y$   is equivalent to conjugating $\phi(0)$ by  an element of $SL(2,\C)$ which yields the upper triangular   entry to be 
                          $$v^2b'+2uv a'-u^2c'$$
                          where $(u,v)$ are homogeneous coordinates in $\PP(V^*_y)$. This quadratic term gives two values where the entry is zero and hence the    Hecke curve intersects the critical locus in two points (corresponding to the two eigenspaces  of $\phi(0)$).
                            
                           The original Higgs field $\Phi$ had just a double zero at $y$ and so is given by the direct image of a line bundle $L$ on the $K(-y)$ spectral curve $\tilde S$.  The fibre $V_y=L_p \oplus L_q$ so if $\alpha\in V^*_y$ annihilates $L_q$, then the bundle $V'$ is defined in terms of $p$ alone and $V'$ is the direct image of the Hecke transform of $L$ at the point $p$. From the exact sequence 
                           $$0\rightarrow {\mathcal O}_{\tilde S}(L(-p))\rightarrow {\mathcal O}_{\tilde S}(L)\rightarrow {\mathcal O}_p(L)\rightarrow 0$$
                            this  is $L(-p)$. Similarly at the other point we get $L(-q)$. It follows that  the two points in the critical locus are related by translation by ${\mathcal O}(p-q)$. Note that since the involution on $\tilde S$ interchanges $p$ and $q$, this lies in the Prym variety.
                           
                           In \cite{Got} the compactification by   generalized parabolic bundles is the smooth projective bundle $\PP(L_p\oplus L_q)$ over the Prym variety of $\tilde S$ which has two sections, the zero  and $\infty$-section. Then the singular fibre is obtained by identifying points in one section with their translates by ${\mathcal O}(p-q)$ on the other section. The critical locus is the result of this identification. It appears  as a normal crossing divisor in the singular fibre.
 (This identification of $0$ and $\infty$ over different fibres was overlooked on page 113 in the author's early paper \cite{Hit2} and is clarified in \cite{HwRam} and \cite{Got}.)

Looking back at the classical case of geodesics on the ellipsoid in Section \ref{classic}, we can view the singular fibre as a real version of the above. The geodesics through $u$ or $v$ are  circle bundles over a circle (instead of a $\PP^1$-bundle over a torus) and the identification by $a-b$ is  provided by the intersection of these two tori in the two geodesics joining $u$ and $v$.

\section{Geometry and topology}
In this section we derive some topological and differential-geometric properties of critical loci using the techniques above.
  \subsection{Multiplicities for the nilpotent cone}
   
   The {\it nilpotent cone} is the subvariety of ${\mathcal M}$ of pairs $(V,\Phi)$ where $\tr\Phi^2\equiv 0$, the fibre of $F$ over $0\in H^0(\Sigma, K^2)$. It has several irreducible components. One is the case $\Phi=0$, the moduli space ${\mathcal N}$ of stable bundles. Here $\tr \Phi\dot\Phi$ is identically zero so this is part of the critical locus. In the other cases $\Phi$ is nilpotent but not identically zero.

          The kernel of a nilpotent $\Phi$ defines an invariant  line bundle $L\subset V$ which by stability has degree $-k<0$ and so the  bundle $V$ is an extension
          $$0\rightarrow L\rightarrow V\rightarrow L^{-1}\rightarrow 0$$
          and a section 
         $b\in H^0(\Sigma, \Hom(L^{-1},LK))=H^0(\Sigma, L^2K)$ is the nilpotent Higgs field. For stability $b\ne 0$  and so $\deg L^2K=2g-2-2k\ge 0$.          
            Different integers $2n=2g-2-2k$ give different components. In terms of Higgs fields vanishing on a  divisor $D$ here we have  ${\mathcal O}(D)\cong {\mathcal O}(L^2K)$.
            
          As in \cite{Hit1}, each component is the closure of a smooth Lagrangian submanifold which is an unramified covering of a vector bundle over the symmetric product $\Sigma^{(2n)}$ of the curve. The zeros  of $b$ define the point in the symmetric product (so that $0\le n< g-1$) and the extension class in $H^1(\Sigma, L^2)$ defining $V$ is the fibre. The covering is the choice of $L$ from the $L^2$ determined by $b$. 
          
         In the  case  $n=0$ $L\cong K^{-1/2}$ and there are $2^{2g}$ components corresponding to the different square roots $K^{1/2}$ of $K$. Then $b\in H^0(\Sigma, L^2K)=H^0(\Sigma, {\mathcal O})$ is non-vanishing so these are not critical points for $F$. 
         The other cases are part of the critical locus and the kernel sheaf of $\ad \Phi$ is clearly $L^{2}\cong K^{-1}(D)$ and the cokernel sheaf $L^{-2}K \cong K^2(-D)$.
         However, since $\Phi$ is everywhere nilpotent, clearly it does not correspond to a nondegenerate critical point. 
         
        Instead, we can use our formula for the Hessian to determine topological information.
   Properness of $F$ means that the closure of the different algebraic components  define compactly supported cohomology classes. Their sum  with appropriate multiplicities is cohomologous to a smooth fibre and these multiplicities are what we can calculate.

          As in the general case, if $a\in H^1(\Sigma,K^{-1})$ is represented by $\bar\partial\alpha/s$  we have a critical point where now $s=b$.  If $b$ has $n$ simple zeros $x_1,\dots, x_n$ then we have $n$ functions $f_i$, for each of which we have a critical point, and from (\ref{hess}) the Hessian of $f_i$ is a multiple of the symmetric bilinear  form $q_k=\tr \dot\Phi_1\dot\Phi_2(x_k)$. 
          
          \begin{prp} For $n\ne 0$  the Hessians of $f_i$ are linearly independent at  a generic point on a component .
          \end{prp}
                \begin{prf} For $g=2$ the inequality $0\le n< g-1$ means $n=0$ and the components are not critical. The statement is an open condition so for $g>2$  assume $\Sigma$ is non-hyperelliptic and hence the canonical map embeds the curve in $\PP^{g-1}$. We shall say that a finite subset $X\subset \Sigma$ is in general position if each subset $Y\subseteq X$ of $k\le g-1$ points spans a $k-1$-dimensional projective subspace which contains no more points of $X$. 
                
                For $n=2g-2-2d$ let $D$ be  an effective divisor consisting of a set of $n$ points in $\Sigma$ in general position. Choose a line bundle $L$ such that $L^2K\cong{\mathcal O}(D)$, and take $V$ to be the trivial extension $V=L\oplus L^*$.  
                
               For each point $x_i\in D$ separate the points  of  $D\backslash  \{x_i\}$  into two disjoint subsets   $D_1,D_2$ with  $\deg D_1=g-1-d, \deg D_2=g-1-d-1$. By general position there exist sections $a_1,a_2$ of $K$ such that $a_1, a_2$ vanish on $D_1, D_2$ but $a_1(x_i)$ and $a_2(x_i)$ are non-zero.  The sections 
               $$\dot\Phi_1=\pmatrix{a_1 & 0\cr
               0 & -a_1}   \qquad \dot\Phi_2=\pmatrix{a_2 & 0\cr
               0 & -a_2}$$
    of $H^0(\Sigma, \End_0 V\otimes K)$ now satisfy the condition   $\tr \dot\Phi_1\dot\Phi_2(x_j)=0$ for $j\ne i$ and     $\tr \dot\Phi_1\dot\Phi_2(x_i)\ne 0.$  
    
    So if $\sum_k c_kq_k=0$ then from this example $c_i=0$ and similarly all $c_k$.

               \end{prf}

                \begin{prp} The multiplicity of a component in the nilpotent cone for which $\Phi$ vanishes on a degree $n$ divisor is $2^n$. The multiplicity of the component ${\mathcal N}$ is $2^{3g-3}.$
                 \end{prp}
                 \begin{prf}   We take a generic point as above. Firstly we know that this is critical for $f_1,\dots, f_n$. The other derivatives $df_{n+1},\dots, df_{3g-3}$ are linearly independent. This follows from 
                 the general observation that since $\tr a\Phi\dot \Phi\in H^1(\Sigma,K)=0$ for all $(\dot A,\dot\Phi)$ at a critical point, it does so in particular for $(\dot A,\dot\Phi)=(0,\dot\Phi)$ where $\dot\Phi$ is any holomorphic Higgs field. It follows that  $a\Phi=0\in H^1(\Sigma, \End_0V)$. In our case 
                 the condition $a\Phi=0$ implies $ab=0\in H^1(\Sigma, L^2)$ since for a direct sum $H^1(\Sigma, L^2)\subset H^1(\Sigma, \End_0V)$. But from the long exact sequence 
                 $$\rightarrow H^0(D, L^2)\rightarrow H^1(\Sigma, K^{-1})\stackrel{b}\rightarrow H^1(\Sigma, L^2)\rightarrow$$
                 the class $a$ is a linear combination of $f_1,\dots, f_n$.
                 
                 We can therefore take a transverse slice $f_i=const$, $n+1\le i\le 3g-3$ and consider the multiplicity of a single point in $\C^n$ where the Hessians are linearly independent. This is the same as the degree of the projective tangent cone which is here the intersection   of $n$ quadrics and by B\'ezout's theorem this is $2^n$, giving the result. 
                 
                 As for the moduli space ${\mathcal N}\subset {\mathcal M}$ where $\Phi\equiv 0$, the cotangent bundle $T^*{\mathcal N}$ is a Zariski open set in ${\mathcal M}$. Take a stable bundle $V$ and 
                 let $a_1,\dots, a_{3g-3}$ be a basis for $H^1(\Sigma, K^{-1})$. Then if the  corresponding Hessians are linearly dependent at $[V]$, there exists $a$ such that $a\tr\Phi_1\Phi_2=0$ for all $\Phi_1,\Phi_2\in H^0(\Sigma,\End_0V\otimes K)$.  However the map $F$ restricted to the Lagrangian cotangent fibre over $[V]\in {\mathcal N}$ in ${\mathcal M}$ surjects onto $H^0(\Sigma, K^2)$ \cite{ Hit1} and  $aq=0$ for all $q$ means $a=0$.

                  Hence the multiplicity is $2^{3g-3}$.

                 \end{prf}

\subsection{The cohomology class of ${\mathcal C}_1$}

 Each critical locus has a closure   $\bar {\mathcal C}_d$ which is a noncompact  analytic cycle,  and 
          is Poincar\'e dual to a  cohomology class in $H^{4d}({\mathcal M},\Z)$.  The cycle $\bar {\mathcal C}_1$ consists of an open set ${\mathcal C}^0_1$ of nonsingular points together with lower-dimensional pieces and as in \cite{GH} the holomorphicity means that integration of a form with compact support over the open piece defines the cohomology class.
          
           The first critical locus is the degeneracy set of the $m$ sections $X_1,\dots,X_m$  of the rank $2m$ tangent bundle $T{\mathcal M}$ where  $m=3g-3$. Cohomology classes of degeneracy sets are in general linked to Chern classes but here we have   Hamiltonian vector fields, which are not generic in the usual algebraic geometric sense. We shall however prove:
            \begin{prp} \label{chern1}The closure $\bar {\mathcal C}_1$ of the first critical locus is Poincar\'e dual to the Chern class $c_2(T{\mathcal M})$.
          \end{prp} 
\begin{prf} 

There are different approaches to this question. 

\noindent 1. One is to follow \cite{Sawon} where the author works in the compact case, using exact sequences of sheaves. In that case the base is $\PP^m$ rather than $\C^m$ which gives a non-zero contribution if $m>1$.

\noindent 2. A second method  is to adapt the approach of  \cite{Bol} where the authors work with a {\it real} integrable system and introduce an almost complex structure compatible with the symplectic structure. The Hamiltonian vector fields then define sections $Y_j=X_j+iIX_j$ of a complex vector bundle and linear dependence is then the vanishing of the complex exterior product $Y_1\wedge\dots\wedge Y_m$. This is  a section of the line bundle $\Lambda^m T^{1,0}$, whose {\it first} Chern class is $c_1(T)$.  The adaptation is to use the hyperk\"ahler structure (though we don't need integrability). 

  The moduli space ${\mathcal M}$ is a hyperk\"ahler manifold which gives a quaternionic structure to $T{\mathcal M}$: an antilinear endomorphism $J$ such that $J^2=-1$.  Since our complex vector fields $X_i$  are commuting Hamiltonian vector fields we have  $\omega(X_i,X_j)=0$. But $\omega(X,JY)$ is the K\"ahler form for complex structure $I$, which means that $JX_i$ is metrically orthogonal to $X_1,\dots, X_m$. 
So the complex vector fields $X_i$ are linearly dependent  if and only if there is a linear dependency over the quaternions  in the $m$-dimensional quaternionic vector bundle $T{\mathcal M}$. 

If $X_1$ vanishes at a point and $X_2,\dots,X_m$ are linearly independent there then we can extend them to a local basis $e_1,\dots,e_m.$ The  transversality condition to make  degeneracy a submanifold of real codimension $4$ is the surjectivity of the projection of the derivative $DX_1$ onto the (quaternionic) coefficient of $e_1$. In our case nondegeneracy provides a  standard model which is the critical locus of the function $z_1z_2$ on $\C^2\times \C^{2m-2}$. Then 
$$X_1=z_1\frac{\partial}{\partial z_1}-z_2\frac{\partial}{\partial z_2}= (z_1-z_2 j)\frac{\partial}{\partial z_1}$$
and $(z_1,z_2)\mapsto z_1-z_2 j$ is an isomorphism.

At this point we have a general situation -- a trivialization of  a principal $G$-bundle $P$ for a simple Lie group $G=Sp(m)$ in the complement of a codimension $4$ submanifold $C\subset M$. Around each point in $C$ we also have a local trivialization in a neighbourhood $U_{\alpha}\cong U\times \R^4\backslash {0}$. This covering by open sets with trivializations gives transition functions  $g_{0\alpha}:(M\backslash C)\cap U_{\alpha}\rightarrow G$. 

If $\Z\cong \pi_3(U_{\alpha})\rightarrow \pi_3(G)\cong \Z$ is an isomorphism then $C$ is Poincar\'e dual to the characteristic class of $P$ corresponding to the generator of $H^4(BG)$. Now $S^3\cong Sp(1)\subset Sp(m)$ generates $\pi_3(Sp(m))$ and 
$(z_1,z_2)\mapsto z_1-z_2 j$ is a diffeomorphism so for ${\mathcal C}_1$ this is a generator. Since $Sp(m)\subset SU(2m)$ is an isomorphism on $\pi_3$ and $c_2$ generates $H^4(BSU(2m))$ we have the result up to a choice of sign which can be fixed by a standard calculation.

\begin{rmk} The above method holds for a general holomorphic integrable system $F:M\rightarrow B$ since $F^*(T^*B)\subset T^*M$ is  generically maximally isotropic. Choosing an almost  quaternionic structure $J$ we have a quaternionic vector bundle homomorphism $F^*(T^*B)\oplus JF^*(T^*B)\rightarrow T^*M$ and, if the critical locus $C\subset M$ is nondegenerate, the local contribution is the same as the above. If we use $\ch_2=(c_1^2-2c_2)/2$ which is additive over direct sums then 
$$[C]=F^*\ch_2(T^*B)+F^*(\ch_2(\bar T^*B))-\ch_2T^*M=c_2(M)-2F^*c_2(B).$$
For example, an elliptic K3 surface has base $\PP^1$ and $c_2(\PP^1)=0$ for dimensional reasons, so a nondegenerate $C$ consists of $c_2(M)=24$ points. 
\end{rmk}

\noindent 3. The final method, in the case of Higgs bundles, is to use the identification of ${\mathcal C}_1$ in terms of Higgs fields which vanish at a point. Over ${\mathcal M}\times \Sigma$ we have  a universal Lie algebra bundle $\End_0V$ and  also a universal Higgs field $\Psi$ which is a section of the rank $3$ bundle $\End_0V\otimes q^*K$ where $q$ is the projection onto the second factor.  It vanishes on a subvariety $Z$ whose projection $p(Z)\subset {\mathcal M}$ is $\bar{\mathcal C}_1$. 

If $Z\subset {\mathcal M}\times \Sigma$ has the expected codimension $3$ it is Poincar\'e dual to $c_3(\End_0V\otimes q^*K)$. For $d=1$, the degenerate cases, such as $\Phi$ having more than one zero, or vanishing entirely, have higher codimension  which means that  $\bar{\mathcal C}_1$ is indeed Poincar\'e dual to $p_*c_3(\End_0V\otimes q^*K)$.

Recall some basic information about the Chern classes of ${\mathcal M}$  \cite{HT}.
     Representing a tangent vector using the first hypercohomology spectral sequence for ${\mathbb H}^1$ and setting $\lie{g}=\End_0V$ we have 
    $$0\rightarrow H^0(\Sigma,\lie{g})\rightarrow H^0(\Sigma,\lie{g}\otimes K)\rightarrow T_x{\mathcal M}\rightarrow H^1(\Sigma,\lie{g})\rightarrow H^1(\Sigma,\lie{g}\otimes K)\rightarrow 0.$$
And so  over ${\mathcal M}\times \Sigma$ the Chern character is 
    $\ch(T{\mathcal M})=\ch(p_{!}(\lie{g}\otimes q^*K-\lie{g}))$.
     By Grothendieck-Riemann-Roch we have 
    $\ch(T{\mathcal M})=(2g-2)p_*(\ch(\lie{g})q^*\sigma)$
    where the class $\sigma\in H^2(\Sigma,\Z)$ is the positive generator.
    
     The rank $3$ vector bundle $\lie{g}$ has second Chern class 
   \begin{equation}
   c_2(\lie{g})=2\alpha\otimes\sigma -\beta\otimes 1+4\sum_{i=1}^g\psi_i\otimes e_i
   \label{c2eq}
   \end{equation}
    where $\beta\in H^4({\mathcal M},\Z)$ and $e_i\in H^1(\Sigma,\Z)$ so that $\sigma c_2(\lie{g})=-\sigma\beta$. This gives 
    $$\ch(T)=(2g-2)\left(3+2\frac{\beta}{2!}+2\frac{\beta^2}{4!}+\dots\right)$$
    which yields the total Chern class 
   $ c(T)=(1-\beta)^{2g-2}$ and $c_2(T)=-(2g-2)\beta.$

    But from (\ref{c2eq}) we obtain  $$c_3(\lie{g}\otimes q^*K)=(2g-2)c_2({\lie g})\sigma=-(2g-2)\beta\sigma$$
    and hence the cohomology class of  ${\bar\mathcal C}_1$ is 
    $$p_*(-(2g-2)\beta\sigma)=-(2g-2)\beta=c_2(T).$$

\end{prf}

\subsection{The normal bundle of ${\mathcal C}^0_1$} \label{normal}

Since ${\mathcal C}^0_1$ is a holomorphic symplectic submanifold of ${\mathcal M}$, its normal bundle $N$ is the symplectic orthogonal of the tangent bundle. It is therefore a holomorphic direct summand in $T{\mathcal M}$ restricted to ${\mathcal C}^0_1$. The symplectic form also makes $N$ an $SL(2,\C)$-bundle and by the result of the previous section $c_2(N)=c_2(T{\mathcal M})$. 

The normal space at a point on the critical locus has appeared implicitly several times in the paper so far and here we put these different viewpoints together.

\noindent 1. Return to the description of the tangent space as the hypercohomology space ${\mathbb H}^1$ and the exact sequence $$0\rightarrow H^1(\Sigma, \ker \Phi)\rightarrow {\mathbb H}^1\rightarrow H^0(\Sigma, \coker \Phi)\rightarrow 0.$$
In Section \ref{Cint} by considering the $K(-D)$-twisted Higgs field $\phi=\Phi/s$ we saw that $\ker\Phi\cong K^{-1}(D)$. The cokernel sheaf $\coker \phi$ was $K^2(-2D)$ so for the Higgs bundle we have an exact sequence 
$$0\rightarrow {\mathcal O}(K^2(-2D))\rightarrow \coker \Phi \rightarrow {\mathcal S}\rightarrow 0$$
for a quotient sheaf ${\mathcal S}$ supported on $D$.

Around a nondegenerate point $\Phi=z\phi$ where $\phi$ is semisimple and then the image of ${\mathcal O}(K^2(-2D))$ by its product with $s$ gives a subsheaf ${\mathcal O}(K^2(-D))\subset {\mathcal S}$ and its quotient ${\mathcal Q}$ is isomorphic to the image of $\ad \phi(0)$. 
If $L$ is the line bundle defining $\phi$ this is the direct sum $L_pL_q^{-1}\oplus L_p^{-1}L_q$. Writing  the tangent space of ${\mathcal C}_d^0$ as
$$0\rightarrow H^1(\Sigma, K^{-1}(D))\rightarrow T{\mathcal C}_d^0\rightarrow H^0(\Sigma, K^2(-D))\rightarrow 0$$
 the normal bundle at this point is 
$$H^0(\Sigma, {\mathcal Q})\cong \bigoplus_{i=1}^d (L_{p_i}L_{q_i}^{-1}\oplus L_{p_i}^{-1}L_{q_i}).$$

\noindent 2. In the proof of Proposition \ref{non} we  identified the subspace $E^{\perp}/E$ for the degeneracy criterion as a direct sum 
$$E^{\perp}/E\cong \bigoplus_{i=1}^d T_i$$
of tangent spaces of coadjoint orbits of $\phi(x_i)$. This is clearly  another view of the above description of the normal bundle.

\noindent 3. In the compactification of a singular fibre by generalized parabolic bundles the smooth $\PP^1$-bundle is immersed in ${\mathcal M}$ with self-intersection the critical locus. The normal bundle is then the direct sum of the normal bundle of the zero section and the infinity section in $\PP(L_p\oplus L_q)$, which is  $L_pL_q^{-1}\oplus L_p^{-1}L_q$.

\noindent 4. The Hecke curves of Section \ref{full} meet, as we observed, a smooth fibre of the integrable system for ${\mathcal C}^0_1$ in two points corresponding to line bundles $L, L(p-q)$. It then also follows that through each point of ${\mathcal C}^0_1$ there pass two Hecke curves. The tangent spaces of these curves span the fibre of the normal bundle at that point.  We can see this more concretely by referring to Equation (\ref{Higgsprime}). 
$$\Phi=\pmatrix{a'z & b'\cr
                            c'z^2 & -a'z}$$
                            which we can write as $\Phi=Z+zB+z^2C$ where $Z$ is nilpotent  and $B$ in its Borel subalgebra. 
                            
                            Hecke transforms are local modifications so a tangent vector  $(\dot A,\dot\Phi)$ along a Hecke curve must have a localized description which is as follows. 
                            Take a nilpotent $Z$ in $H^0(D,\End_0 V'(D))$ corresponding to the distinguished subspace of $V'$. Then, as in \cite{NRam1}, the connecting homomorphism $\delta: H^0(D,\End_0 V'(D))\rightarrow H^1(\Sigma, \End_0V')$ defines the infinitesimal deformation $\delta(Z)$ of the holomorphic structure along a tangent to the Hecke curve. A Dolbeault representative $\dot A\in \Omega^{01}(\Sigma, \End_0 V')$ is defined 
                             in the usual way by an element   
                            $\dot A=\bar\partial  \beta/s$                         
where $\beta$ is supported near $D$ and restricted to $D$ is $Z$. Then locally 
$$[\dot A,\Phi] = \frac{1}{z}[\bar\partial  \beta,Z+zB+z^2C]=\bar\partial ([\beta, B]+z[\beta, C])$$
since $\beta$ commutes with $Z$ at $z=0$. 
Putting $\dot\Phi$ equal to the smooth section $-[\beta,\Phi]$ we have $\bar\partial\dot\Phi+[\dot A,\Phi]=0$ and a tangent vector to ${\mathcal M}$. 

The semisimple $\phi(0)$ determines two nilpotents $Z^+, Z^-$ spanning the image of $\ad\phi$ and hence two normal directions which span the normal bundle at this point.

\subsection{The induced metric on ${\mathcal C}_1$}

The moduli space ${\mathcal M}$ is not only a complex symplectic manifold, but also has a natural hyperk\"ahler metric as a consequence of the existence of solutions of the Higgs bundle equations \cite{Hit1}. The critical loci ${\mathcal C}_d^0$ are complex symplectic submanifolds so it is natural to ask whether they are hyperk\"ahler submanifolds. These would then be the support of BBB-branes \cite{KW}.

For a submanifold to be  hyperk\"ahler  is a strong condition -- it means it is holomorphic with respect to all complex structures $I,J,K$ of the hyperk\"ahler family. It is well-known that this means the submanifold is totally geodesic. In fact, the second fundamental form $S(X,Y)$ of a complex submanifold of a K\"ahler manifold is a symmetric complex linear form with values in the normal bundle. So for a hyperk\"ahler submanifold 
$$S(IX,JY)=IJS(X,Y)=JIS(X,Y)=-IJS(X,Y)$$
and so $S$ vanishes.

We shall show that   ${\mathcal C}_1$ is not a hyperk\"ahler submanifold of ${\mathcal M}$ in the case where, for simplicity, we take  the degree of $V$ to be  odd and $g>2$.

Since the hyperk\"ahler metric is unknown we need to prove this by contradiction using some holomorphic consequences of the assumption that it is hyperk\"ahler. If a complex submanifold of a K\"ahler manifold is totally geodesic, then   the normal bundle is a holomorphic direct summand of $TM$ restricted to the submanifold. We have seen in Section \ref{normal} that this  in fact holds in our situation but we can still follow the same idea. 

The involution $(A,\Phi)\mapsto (A,-\Phi)$ is an isometry on the Higgs bundle moduli space and its fixed point set is therefore totally geodesic. If ${\mathcal C}^0_1$ is totally geodesic then its intersection with this fixed point set is again totally geodesic. The fixed point set has different components which  are described in \cite{Hit1} and correspond to representations of the fundamental group of $\Sigma$ into a real form, but we shall focus on the component described by 
$$V=L\oplus L^*\Lambda^2V\qquad \Phi= \pmatrix { 0 & b\cr
                                                                                    c & 0}$$
where $\Lambda^2V$ is a fixed line bundle of degree $1$, $L$ is a line bundle of degree $(g-1)$, $c$ a holomorphic section of $L^{-2}K\otimes \Lambda^2V$ and $b$ of $L^2K\otimes \Lambda^2V^*$. By stability $c$ is nonzero and $L^{-2}K\otimes \Lambda^2V$ is of degree $1$ so $c$ vanishes at a unique point $x\in \Sigma$.  Then $b$ lies in the $3g-4$-dimensional space $H^0(\Sigma, K^2(-x))$ since $\det \Phi=-bc$. 
The intersection with ${\mathcal C}^0_1$ is where $b$ has a simple zero at $x$.

This component is a smooth submanifold  ${\mathcal M}_1\subset {\mathcal M}$ and  (see \cite{Hit1}) the projection to the zero of $c$ gives it the structure of a $2^{2g}$-fold covering of a vector bundle $W$ over $\Sigma$ (the covering is the choice of $L$ given that $L^{-2}K\otimes \Lambda^2V\cong {\mathcal O}(x)$.)

There are points in this submanifold (if $b$ has a multiple zero at $x$) which are not nondegenerate critical points but ${\mathcal M}_1$ is still smooth and moreover if the open set ${\mathcal C}_1^0\cap {\mathcal M}_1$ is totally geodesic so will be the smooth points of its closure. 

So far we have replaced the more complicated geometry of ${\mathcal M}$ by that of a vector bundle over $\Sigma$ but this is still non-compact. We now use the isometric circle action $(A,\Phi)\mapsto (A, e^{i\theta}\Phi)$ which preserves both the critical loci and ${\mathcal M_1}$. Removing the zero section  of $W\rightarrow \Sigma$ it acts freely by scalar multiplication on each fibre $H^0(\Sigma, K^2(-x))$. Its K\"ahler quotient is the $\C^*$-quotient which is the projective bundle $\PP(W)\rightarrow \Sigma$. We now need the following:

\begin{lemma} \label{quotlemma} Let $M$ be a K\"ahler manifold with a free circle action preserving the K\"ahler structure. If $N\subset M$ is an invariant complex totally geodesic submanifold then the K\"ahler quotient of $N$ is totally geodesic in the K\"ahler quotient of $M$. 
\end{lemma}
\begin{prf} Recall that if $\mu$ is the moment map for the action then the  K\"ahler quotient is $\mu^{-1}(0)/S^1$ with its induced metric.

Let $X$ be the vector field generated by the circle action on $M$. Then the normal to the zero set $\mu^{-1}(0)$ of the moment map is $IX$. But $N\subset M$ is $S^1$-invariant and complex so at a point on $N$ the vector field  $IX$ is tangential to $N$. It is therefore normal to the zero set of $\mu$ restricted to $N$.

The Levi-Civita connection for the induced metric on $\mu^{-1}(0)$ is given by removing the normal component:
$$\nabla_YZ-\frac{1}{(X,X)}(\nabla_YZ, IX)IX.$$
But if $N\subset M$ is totally geodesic, $\nabla_YZ$ for vector fields $Y,Z$ on $N$ is the Levi-Civita connection on $N$. Since $IX$ is the normal in $N$ to $N\cap \mu^{-1}(0)$ the above expression gives the Levi-Civita connection on $N\cap \mu^{-1}(0)$ so this is totally geodesic in  $\mu^{-1}(0)$.

Geodesics on the quotient $\mu^{-1}(0)/S^1$ lift to horizontal geodesics $((\dot\gamma, X)=0$) in $\mu^{-1}(0)$. Since $N\cap \mu^{-1}(0)$ is  totally geodesic in $\mu^{-1}(0)$, horizontal geodesics in the submanifold $N\cap \mu^{-1}(0)$ are horizontal geodesics in $\mu^{-1}(0)$.
\end{prf}
\begin{rmk} An analogue of the first part of the proof is the classical situation of a surface $S$ in $\R^3$. Its intersection with a (totally geodesic) plane is not in general a geodesic but if the normal to the surface lies in the plane along the curve of intersection, then the normal to the plane curve is normal to the surface and hence is a geodesic.
\end{rmk}

The $2^{2g}$-fold covering consists of the action $V\mapsto V\otimes U$ for a line bundle $U$ of order $2$. This is a free  isometric action and so the critical locus in ${\mathcal M}_1$ will be totally geodesic if and only if the corresponding subspace of $\PP(W)$ is. We shall investigate this next.

Consider the trivial vector bundle $\Sigma\times H^0(\Sigma,K^2)$ and the evaluation map to $K^2$. The kernel over a point $x\in \Sigma$ consists of sections which vanish at $x$ i.e. $b\in H^0(\Sigma, K^2(-x))$ so the vector bundle $W$ above in the description of ${\mathcal M}_1$ is the kernel. If $g>2$ then the one-jet of a section gives a surjective homomorphism (the bicanonical system is an embedding) from $\Sigma\times H^0(\Sigma,K^2)$ to $J^1(K^2)$ and the kernel $(W_1)_x$ is the space of $b\in H^0(\Sigma, K^2(-x))$ which vanish at $x$. The bundle of $1$-jets is an extension
$$0\rightarrow K^3\rightarrow J^1(K^2)\rightarrow K^2\rightarrow 0$$
where the projection to $K^2$ is the $0$-jet, so $W_1\subset W$ with quotient $K^3$.

 By Lemma \ref{quotlemma}, if ${\mathcal C}_1^0$ is a hyperk\"ahler submanifold of ${\mathcal M}$ then $\PP(W_1)$ is a totally geodesic submanifold of $\PP(W)$.
\begin{prp} The normal bundle of $\PP(W_1)$ in $\PP(W)$ is not a holomorphic direct summand. 
\end{prp}
\begin{prf} Let $H$ be the hyperplane line bundle over  $\PP(W)$. Since $\PP(W_1)\subset \PP(W)$ is an embedding of projective bundles the normal bundle is the normal bundle along the fibres and since $W/W_1\cong K^3$ it is $Hp^*K^3$ where $p:\PP(W_1)\rightarrow \Sigma$ is the projection.

If this is a direct summand then $Hp^*K^3$ injects into the tangent bundle $T\PP(W)$ restricted to $\PP(W_1)$. The projection to $p^*T\Sigma$ is zero since $H^{-1}$ has no sections on the projective space fibres. It therefore must map to the tangent bundle along the fibres $T_F$ for which we have the Euler sequence 
$$0\rightarrow {\mathcal O}\rightarrow p^*W\otimes H \rightarrow T_F\rightarrow 0.$$
Now $H^1(\PP(W_1), H^{-1}p^*K^{-3})=0$ since $H^q(\PP^n,{\mathcal O}(-1))=0$ for all $q$ so the homomorphism lifts to $p^*W\otimes H$. But 
$H^0(\PP(W_1), p^*(W\otimes K^{-3}))\cong H^0(\Sigma, W\otimes K^{-3})$
and  $W$ is a subbundle of the trivial bundle. Since $K^{-3}$ has no holomorphic sections   $H^0(\Sigma, W\otimes K^{-3})=0$.
\end{prf}

We deduce that ${\mathcal C}_1^0$ is not a hyperk\"ahler submanifold. 
 \subsection{The extreme case $d=2g-2$}\label{extreme}
Consider now the induced metric at the other extreme from ${\mathcal C}_1$ -- where  the divisor  $D$ has degree $2g-2$. Since $K^2$ has degree $4g-4$, a quadratic differential with $2g-2$ double zeros must be of the form $q=s^2$ where $s$ is either a section of $K$ or $KU$ where $U^2$ is trivial. In the first case, as we have seen before, $\phi\in H^0(\Sigma, \End_0V)$ contradicts stability, so $U$ must be non-trivial. Then the  critical locus ${\mathcal C}_{2g-2}$   has $2^{2g}$ different components for the choice of $U\in H^1(\Sigma,\Z_2)$ and each is an integrable system of dimension $6g-6-2(2g-2)=2g-2$. The quadratic differentials in the base are of the form $q=\lambda s^2$ and $\dim H^0(\Sigma, KU)=g-1$. 

The spectral curve $\tilde S$ of $\phi=\Phi/s$ is then the unramified covering of $\Sigma$ defined by the class of $U$ in $H^1(\Sigma,\Z_2)$.

\begin{prp}  The component of the critical locus ${\mathcal C}_{2g-2}$ corresponding to $U\in H^1(\Sigma,\Z_2)$ is the fixed point set of the action  $(V,\Phi)\mapsto (V\otimes U,\Phi)$ on the moduli space of Higgs bundles ${\mathcal M}$. It is a flat hyperk\"ahler submanifold of ${\mathcal M}$.

\end{prp}
\begin{prf} 

\noindent 1. We have $V=\pi_*L$ for $L$ a line bundle on $\tilde S$. But $V\otimes U=\pi_*(L\pi^*U)$ and $\pi^*U$ is trivial on $\tilde S$ since $U$ defines the double covering. Hence $V\otimes U\cong \pi_*L=V$, and $\tilde x$, hence $\Phi$, is unchanged.

\noindent 2. The operation $V\mapsto V\otimes U$ is holomorphic in the Higgs bundle complex structure but also, since $U$ has a flat connection with holonomy in the centre $\pm 1$ of $SU(2)$, in the complex structures $J$ and $K$. The complex structure $J$ is   determined by solving the Higgs bundle equations $F+[\Phi,\Phi^*]=0$ and taking the flat $SL(2,\C)$-connection $\nabla_A+\Phi+\Phi^*$.  The  fixed point set is therefore holomorphic with respect to all complex structures and so the induced metric is hyperk\"ahler.

\noindent 3. For an unramified covering the direct image takes a flat connection to a flat connection and a metric to a metric. It follows that the solution to the Higgs bundle equations for $(V,\Phi)$ is the direct image of a solution in the line bundle case on the covering $\tilde S$. The natural metric here is the flat metric on the cotangent bundle of the Jacobian of $\tilde S$ and its restriction to the anti-invariant part, the cotangent bundle of the Prym variety, is flat too. 

\end{prf}

\section{Genus 2}
\subsection{The integrable system}
In the case of a curve $\Sigma$ of genus $2$, there are natural coordinates which, thanks to the work of \cite{vG} and \cite{Gaw}, enables one to write down the integrable system explicitly. 

The background is the original paper \cite{NRam}. For each stable $SL(2,\C)$-bundle $V$, the line subbundles of degree $-1$ are described by a $2\Theta$-divisor in $\Pic^1(\Sigma)$ and then the moduli space of S-equivalence classes of semistable bundles is isomorphic to $\PP(H^0(\Pic^1(\Sigma), 2\Theta))=\PP^3$. The set of equivalence classes of semistable bundles  is given  by bundles of the form $V=U\oplus U^{-1}$ where $U$ is a line bundle of degree zero. The rank $2$ bundle does not distinguish the order and so the locus is $\Pic^0(\Sigma)/\pm 1$ which is the  Kummer quartic surface. It has 16 singularities at the points where $U\cong U^{-1}$, the line bundles of order $2$.

The cotangent bundle of the stable locus lies in the Higgs bundle moduli space as a dense open set and the cotangent bundle of a projective space $\PP(V)$ can be described as the symplectic quotient of $V\times V^*$ by the $\C^*$-action $(p,q)\mapsto (\lambda p,\lambda^{-1}q)$ so that 
$$T^*\PP(V)=\{(p,q)\in V\times V^*:p\ne 0, \langle p, q\rangle =0\}/\C^*.$$

In these coordinates the function $F:{\mathcal M}\rightarrow H^0(\Sigma, K^2)$ is given in \cite{Gaw} by the beautifully symmetric expression
\begin{equation}
F(p,q)=-\frac{1}{128 \pi^2}\sum_{i\ne j}\frac{\langle \sigma(ij)p,q\rangle^2}{(z-z_i)(z-z_j)}dz^2.
\label{formula}
\end{equation}
Here the hyperelliptic curve $p:\Sigma\rightarrow \PP^1$ is $w^2=f(z)=(z-z_1)\dots(z-z_6)$ 
and $\sigma(ij)$ is an element in the finite Heisenberg group   which acts on $V=H^0(\Pic^1(\Sigma), 2\Theta)$. In terms of vector bundles it is the group $\Gamma\cong \Z_2^4$ of line bundles of order $2$ acting as $V\mapsto V\otimes U$ and a central extension by $\pm 1$ acts on $V$. The Weierstrass points $x_i,x_j\in \Sigma$ giving the branch points $z_i,z_j\in \PP^1$ differ by a divisor class of order $2$ and this is $\sigma(ij)$. The power of $2$ in  $\langle \sigma(ij)p,q\rangle^2$ means that the expression is independent of the lift to the central extension and the dual pairing $\langle p,q\rangle$ is clearly $\C^*$-invariant. 

The canonical bundle of $\Sigma$ is $K\cong p^*{\mathcal O}(1)$ and has a basis of sections $dz/w,zdz/w$. The three-dimensional space of  sections of $K^2$ is generated by these and hence of the form 
$$(a_0+a_1z+a_2z^2)\frac{dz^2}{w^2}=\frac{a_0+a_1z+a_2z^2}{f(z)}dz^2.$$
The partial fractions on the right hand side of (\ref{formula}) reduce to this form taking account of the symmetry of the action.

However elegant the formula (\ref{formula}) may be, it is not at all easy to do any calculations with. A more amenable, but less symmetric, expression can be found in \cite{Lor}.  
The authors take the six distinct branch points to be $z_1,z_2,z_3 = 0, 1, \infty$ and $z_4=r, z_5=s, z_6=t$. If $(u_0,u_1,u_2)$ are affine coordinates on $\PP^3$ and $\eta_0du_0+\eta_1du_1+\eta_2du_2$ the canonical one-form on $T^*\PP^3$ then $F$ maps to  $(h_0+h_1z+h_2z^2)dz^2/f(z)$ given (up to a universal factor) by these equations:
\begin{eqnarray*}
h_0&=& rst[\eta_0(u_0^2-1)+\eta_1(u_0u_1+u_2)+\eta_2(u_2u_0+u_1)]^2-\\
&& st[\eta_0(u_0u_1-u_2)+\eta_1(u^2_1+1)+\eta_2(u_1u_2+u_0)]^2+\\
&&4rs(\eta_0u_0+\eta_1u_1)^2-rt[\eta_0(u_0^2+1)+\eta_0(u_0u_1+u_2)+\eta_2(u_2u_0-u_1)]^2\\[7pt]
h_1&=&t(u_0^2+u_1^2+u_2^2+1)[(\eta_0^2+\eta_1^2+\eta_2^2)+(\eta_0u_0+\eta_1u_1+\eta_2u_2)^2]+\\
&&st(u_0^2-u_1^2+u_2^2-1)[(\eta_0^2-\eta_1^2+\eta_2^2)-(\eta_0u_0+\eta_1u_1+\eta_2u_2)^2]+\\
&&4r(u_0u_2-u_1)[\eta_0\eta_2+(\eta_0u_0+\eta_1u_1+\eta_2u_2)\eta_1]+\\
&&4sr(u_2u_0+u_1)[\eta_2\eta_0-(\eta_0u_0+\eta_1u_1+\eta_2u_2)\eta_1]+\\
&&4s(u_1u_2+u_0)[\eta_1\eta_2-(\eta_0u_0+\eta_1u_1+\eta_2u_2)\eta_0]+\\
&&4rt(u_0u_1+u_2)[\eta_0\eta_1-(\eta_0u_0+\eta_1u_1+\eta_2u_2)\eta_2]\\[7pt]
h_2&=&s[\eta_0(u_2u_0+u_1)+\eta_1(u_1u_2+u_0)+\eta_2(u_2^2-1)]^2-\\
&&[\eta_0(u_2u_0-u_1)+\eta_1(u_1u_2+u_0)+\eta_2(u_2^2+1)]^2-\\
&&t[\eta_0(u_0u_1+u_2)+\eta_2(u_1u_2-u_0)+\eta_1(u_2^2+1)]^2+4r(\eta_1u_1+\eta_2u_2)^2
\end{eqnarray*}
The Kummer surface in these coordinates has the equation
\begin{dmath} 
t(s-1)(u_0^4+u_1^4+u_2^4+1)-8(r(s-t+1)-s)u_0u_1u_2\\
-2(st+t-2s)(u_1^2u_2^2+u_0^2)
- 2(s-1)(2r-t)(u_2^2u_0^2+u_1^2)
+t(2r-(s+1))(u_0^2u_1^2+u_2^2)=0. 
\label{Kummer}
\end{dmath}
We shall try and describe the geometry of the subintegrable systems ${\mathcal C}_1$ and ${\mathcal C}_2$ by the algebraic geometry of $\PP^3$, and in particular the image of the abelian varieties, but also refer to these formulas when appropriate.

\subsection{The critical locus ${\mathcal C}_1$}
Sections of $K^2$ are, as we have seen,  of the form $(a_0+a_1z+a_2z^2)dz^2/{f(z)}$ and a double zero occurs in two ways. 
In case $a_0+a_1z+a_2z^2=a_0(z-c)^2$, we have two double zeros $(w,z)=(\pm f(c),c)$ and then $D$ is  a canonical divisor and  there is no stable critical locus. 

The alternative is that $z=z_1$, say, and $D=x_1$, a ramification point of $p:\Sigma\rightarrow \PP^1$. Now $2x_1\cong p^*{\mathcal O}(1)\cong K$ so $K(-D)\cong K^{1/2}\cong{\mathcal O}(x_1)$, an odd theta characteristic. Then $a_0+a_1z+a_2z^2=(az+b)(z-z_1)$.

Let $s$ be the section of $K^{1/2}$ vanishing at $x_1$, so that $s^2=z-z_1$. Now $w$ is a section of $K^3$ which satisfies $w^2=f(z)=(z-z_1)f_1(z)$ and vanishes with multiplicity $1$ at $x_1$, so $w/s$ is a well-defined section of $K^{5/2}$.

Since $K^2(-2D)\cong K$ the $K(-D)$-twisted  spectral curve $\tilde S$ of genus $4$ has equation
$\tilde x^2=(az+b)$ where $\tilde x$ is a section of $\pi^*K^{1/2}$. Hence $\tilde w=\tilde x w/s$ is a section of $\pi^*K^3$ and 
$$\tilde w^2=(az+b)(z-z_2)\dots (z-z_6)$$
 and so we have a map $\pi_1:\tilde S\rightarrow \Sigma_1$ to another genus $2$ curve $\Sigma_1$.  This is the quotient by the involution $(\tilde x, w)\mapsto (-\tilde x,-w)$ on $\tilde S$. Hence $\tilde S$ has two commuting involutions $\sigma,\sigma_1$ with quotients $\Sigma, \Sigma_1$.
  The projection $\pi:\tilde S\rightarrow \Sigma$ has two branch points $z=-b/a, w=\pm\sqrt{ f(-b/a)}$ and the other  $\pi_1:\tilde S\rightarrow \Sigma_1$ has two branch points $z=z_1, \tilde w=\pm \sqrt{(az_1+b)f'(z_1)}$. 
 
 In this case the Prym variety $\Prym(\tilde S,\Sigma)$ is isomorphic to $\pi_1^*\Pic^0(\Sigma_1)$ so the vector bundle is given by 
 $$V=\pi_*(\pi_1^*L) \otimes K^{1/2}.$$
 and the fibre of the integrable system projects to ${\mathcal N}=\PP^3$ as the Kummer surface for the curve $\Sigma_1$. 
\begin{rmk} It is rare for a Prym variety to be a Jacobian but here the map $\pi_1:\tilde S\rightarrow \Sigma_1$ is the exceptional case of degree $2$ ramified at two points, as in  \cite{Birk}, 12.3.3. 
\end{rmk}
  As $(a,b)\in \PP^1$ varies we have a pencil of singular Kummer quartics. The formula (\ref{Kummer}) in the previous section for the general Kummer surface is linear in $r$ and this is the pencil.
\begin{rmks} 

\noindent 1.  Higgs bundles twisted by a theta-characteristic form the topic of W.Oxbury's 1987 thesis \cite{Ox}. He characterizes the pencil above as the quartics which meet the Kummer surface of $\Sigma$ tangentially in a certain smooth curve. This curve is nonsingular, genus $5$ and degree $8$ and passes through all $16$ singularities. It is  defined by $V=L\oplus L^{-1}$ where $L$ has degree zero and $H^0(\Sigma, L^2K^{1/2})=1$.

\noindent 2.  Given the identification of the Prym variety with $\Pic^0(\Sigma_1)$ we can also ask for the image of the Hecke curves, the rest of the singular locus. From \cite{NRam1}, \cite{Hw} these project in general to rational curves of degree $4$ with respect to the anticanonical bundle of ${\mathcal N}$ so they are lines in $\PP^3$. We have seen in Section \ref{full} that they intersect a smooth fibre in two points which differ by the divisor $p-q$. In our case $p$ and $q$ are the ramification points in $\tilde S$ for the projection $\pi_1$ and so they map to $z=z_1, \tilde w=\pm \sqrt{(az_1+b)f'(z_1)}$ in $\Sigma_1$. Let $U$ be the degree zero line bundle on $\Sigma_1$ given by the difference of these two, then the two Hecke curves through $V_L=\pi_*(\pi_1^*L) \otimes K^{1/2}$ are the lines joining $V_L$ to $V_{LU}$ and $V_{LU^{-1}}$. 
  \end{rmks}

  We conclude that in this case ${\mathcal C_1^0}$ has six components corresponding to the choice of $z_i$ and the family of abelian varieties is described by their Kummer surfaces which lie in the pencil.

\subsection{The critical locus ${\mathcal C}_2$}

This is the case of $d=2g-2$ dealt with in Section \ref{extreme}, and the locus is the fixed point set of an element in $\Gamma=H^1(\Sigma, \Z_2)=\Z_2^4$. This action is described in \cite{Lor} and we may 
consider a typical  involution $(u_0,u_1,u_2,\eta_0,\eta_1,\eta_2)\mapsto (u_0,-u_1,-u_2,\eta_0,-\eta_1,-\eta_2)$.  

Then the fixed point set in $\PP^3$ consists of two lines, corresponding to the two components of the Prym variety of an unramified cover. One of them is on the plane at infinity and the other 
$u_1=u_2=0$. We focus on this one. The fixed point set in ${\mathcal M}$ is $u_1=u_2=\eta_1=\eta_2=0$ and so in the formulae above 
 the Hamiltonian functions are  
\begin{eqnarray*}
h_0&=&\eta_0^2r[st(u_0^2-1)^2
+4su_0^2-t(u_0^2+1)^2]\\
h_1&=&\eta_0^2[t(u_0^2+1)^2
-st(u_0^2-1)^2
-4su_0^2]=h_0/r\\
h_2&=&0
\end{eqnarray*}
So the locus maps to the one dimensional subspace of $H^0(\Sigma, K^2)$ spanned by the quadratic differential $(z-r)dz^2/{f(z)}$ and the fibre is defined by 
$$y^2=st(x^2-1)^2
+4sx^2-t(x^2+1)^2.$$
The discriminant of the quartic is $4096s^2t^2(s-1)^2(s-t)^2(t-1)^2$ and so, as expected, the fibre is a nonsingular  elliptic curve. 
  
\vskip 1cm
 {Mathematical Institute,
Radcliffe Observatory Quarter,
Woodstock Road,
Oxford, OX2 6GG}

\begin{thebibliography}{11}
        
  \bibitem{Alex}
V.Alexeev, 
{\it Compactified Jacobians and the Torelli map}, 
Publ. Res. Inst. Math. Sci. {\bf 40 }(2004) 1241Ð-1265. 
%
\bibitem{BNR}
A.Beauville, M.S.Narasimhan \& S.Ramanan, 
{\it Spectral curves and the generalised theta divisor},
J. Reine Angew. Math. {\bf  398} (1989), 169Ð-179.
%
\bibitem{Birk}
 C.Birkenhake \& H.Lange, ``Complex abelian varieties", Second edition, Grundlehren der Mathematischen Wissenschaften {\bf  302}. Springer-Verlag, Berlin, (2004).
%
\bibitem{Bol}
A.V. Bolsinov \& A.A.Oshemkov,
{\it Singularities of integrable Hamiltonian systems}, 
in ``Topological methods in the theory of integrable systems", 1 -- 67, (eds A. V. Bolsinov et al) Camb. Sci. Publ., Cambridge, (2006).  
%
\bibitem{Bis}
I.Biswas \& S.Ramanan, 
{\it An infinitesimal study of the moduli of Hitchin pairs},
J. London Math. Soc. {\bf 49} (1994)  219 Ð- 231. 
%
\bibitem{Gaw}
K.Gaw\c{e}dzki  \& P.Tran-Ngoc-Bich, 
{\it Self-duality of the $SL_2$ Hitchin integrable system at genus $2$},
Commun.Math.Phys. {\bf 196} (1998)  641--670. 
%
\bibitem{Got}
P.B.Gothen \& A.G.Oliveira (2013) {\it The singular fiber of the Hitchin map,}
Int. Math. Res. Not. {\bf 2013} (2013) 1079Ð-1121.

%
\bibitem{GH}
P.Griffiths \& J.Harris,
``Principles of Algebraic Geometry", John Wiley \& Sons, New York,  (1978).
%
\bibitem{HT}
T.Hausel \& M.Thaddeus,
{\it Relations in the cohomology ring of the moduli space of rank 2 Higgs bundles,} Journal of the AMS, {\bf 16} (2003) 303--327.

%
 \bibitem{Hit1}
 N.J.Hitchin, {\it The self-duality equations on a Riemann surface,} Proc. London
Math. Soc. {\bf 55} (1987),  59--126.
%
\bibitem{Hit2}
 N.J.Hitchin, {\it Stable bundles and integrable systems}, Duke Math. J. {\bf 54} (1987) 91Ð-114.
%
\bibitem{Hit3}
 N.J.Hitchin, {\it Remarks on Nahm's equations,} in ``Modern Geometry: A celebration of the work of Simon Donaldson'', Proc of Symposia in Pure Mathematics, American Math. Soc.  (to appear) arXiv:1708.08812
 %
 \bibitem{Hw}
 J-M Hwang,
{\it Tangent vectors to Hecke curves on the moduli space of rank 2 bundles over an algebraic curve},   Duke Math. J. {\bf 101} (2000) 179--187.

 %
\bibitem{HwRam}
J-M Hwang \& S.Ramanan,
{\it Hecke curves and Hitchin discriminant},  Ann. Sci. \'Ecole Norm. Sup. {\bf  37} (2004) 801--817. 
%
 \bibitem{Ito}
H.Ito, {\it Action-angle coordinates at singularities for analytic integrable systems,} Math. Z. {\bf 206} (1991),  363--407.
%
\bibitem{KW}
A.Kapustin \& E.Witten, {\it Electric-magnetic duality and the geometric Langlands program}, Commun. Number Theory Phys. {\bf 1} (2007) 1--236.
%
\bibitem{Lor}
F.Loray \& V.Heu,
{\it Hitchin Hamiltonians in genus 2,} ``Analytic and Algebraic Geometry" (A.Aryasomayajula et al (eds.)), Springer Nature Singapore Pte Ltd. and Hindustan Book Agency,  (2017) 153--172. arXiv: 1506.02404v1
%
\bibitem{NRam}
M.S.Narasimhan \& S.Ramanan,
{\it Moduli of vector bundles on a compact Riemann surface}, Ann. of Math. {\bf 89} (1969) 14--51. 
%
\bibitem{NRam1}
M.S.Narasimhan \& S.Ramanan,
{\it Geometry of Hecke cycles I}, in ``C.P.Ramanujam -- a tribute" Tata Inst. Fund.
Res. Studies in Math. {\bf 8}, Springer-Verlag, Berlin (1978)  291Ð-345.
%
\bibitem{Nit}
 N.Nitsure, {\it Moduli space of semistable pairs on a curve}, Proc. London Math. Soc. {\bf 62} (1991)  275Ð-300.
%
\bibitem{Ox}
W.Oxbury, {\it Stable bundles and branched coverings over Riemann surfaces}, DPhil thesis, University of Oxford (1987).
%
\bibitem{Sawon}
J.Sawon, {\it On the discriminant locus of a Lagrangian fibration}, Math. Ann. {\bf 341} (2008),  201--221.
%
\bibitem{Sim}
C.T.Simpson, {\it Moduli of representations of the fundamental group of a smooth projective variety II}, Pub. math. IHES {\bf 79}  (1994) 47--129.
%
\bibitem{vG}
B.van Geemen \& E.Previato, {\it On the Hitchin system,} Duke Math.J. {\bf 85} (1996),  659--683.
%
\bibitem{CdV}
Y.Colin de Verdie\`re \& S.Vu Ngoc, 
{\it Singular Bohr-Sommerfeld rules for 2D integrable systems},  
Ann. Sci. \'Ecole Norm. Sup. {\bf  36} (2003) 1Ð-55. 
%
\bibitem{Vey}
J.Vey, {\it  Sur certain syst\`emes dynamiques s\'eparables}  Amer. J. Math.
{\bf 100}  (1978) 591 -- 614.

  \end{thebibliography}
\end{document}